\documentclass[11pt]{amsart}
\usepackage{amsmath,amssymb,mathtools}
\usepackage[margin=1.1in]{geometry}

\usepackage{xcolor}
\usepackage[shortlabels]{enumitem}
\usepackage[colorlinks,linkcolor=magenta,citecolor=blue,urlcolor=blue]{hyperref}

\newtheorem{theorem}{Theorem}[section]
\newtheorem{proposition}[theorem]{Proposition}
\newtheorem{lemma}[theorem]{Lemma}
\newtheorem{corollary}[theorem]{Corollary}
\theoremstyle{definition}
\newtheorem{conjecture}[theorem]{Conjecture}
\numberwithin{equation}{section}

\newcommand{\C}{\mathbb C}
\newcommand{\R}{\mathbb R}
\newcommand{\B}{\mathcal B}
\newcommand{\HH}{\mathcal H}
\newcommand{\KK}{\mathcal K}
\newcommand{\tr}{\operatorname{tr}}

\title[]{A counterexample to the Jordan tuple conjecture of Ball and Helton}
\author{E. Evert and S. McCullough}
\address{Department of Mathematics\\
University of Florida\\
Box 118105\\
Gainesville, FL 32611-8105 }
\email{eevert@ufl.edu, sam@ufl.edu}

\begin{document}
\begin{abstract}
This article constructs a pair of bounded commuting $3$-symmetric operators with no
Jordan-family extension. The operators are coordinate multipliers on a
Sobolev space of $\C^3$-valued functions with an equivalent Hilbert norm.
Their quadratic symbol has no operator polynomial factorization, even with
an auxiliary target Hilbert space. The construction transfers the
nonfactorability of the matrix quadratic used by Russo to the symbol of a
commuting $3$-symmetric tuple. The factorization approach of
McCullough--Russo and Russo then supplies the obstruction to a
Jordan-family extension.
\end{abstract}


\maketitle

\section{The conjecture and the factorization obstruction}
\label{sec:extensions}

For a commuting pair $T=(T_1,T_2)$ in $\B(\HH)$, define its symbol by
\begin{equation}\label{eq:symbol}
 Q_T(s,t)=e^{-itT_2^*}e^{-isT_1^*}e^{isT_1}e^{itT_2},
 \qquad (s,t)\in\R^2.
\end{equation}
The pair is a \emph{commuting $3$-symmetric tuple} if $Q_T$ is a
polynomial of total degree at most two. In the terminology of
\cite{BH80}, such a pair is a \emph{coadjoint family}.

The single-operator theory has its origins in Helton's work
\cite{Helton71,Helton-multx}, which connects Jordan extensions and
multiplication operators on Sobolev spaces with Sturm--Liouville
operators and classical conjugate point theory.  See also \cite{BF90}. It was decisively concluded by Agler
\cite{Agler-thesis}. See also
\cite{Agler-disc}. For commuting families, Ball and Helton
\cite[Part~III]{BH80} relate simultaneous Jordan extensions to
higher-dimensional disconjugacy and conjugate surface theory.

A Jordan family of order two can be written as $J_j=S_j+N_j$ on a Hilbert
space $\KK$, where
\begin{equation}\label{eq:jordan-family}
 S_j=S_j^*,\qquad S_iS_j=S_jS_i,\qquad
 S_iN_j=N_jS_i,\qquad N_iN_j=0
 \quad (i,j\in\{1,2\}).
\end{equation}
These are the Jordan families of Ball--Helton; the commutation relations
between the selfadjoint and nilpotent parts follow from their structure
result for commuting Jordan operators \cite{BH80}.
A Jordan-family extension of $T$ means a family satisfying
\eqref{eq:jordan-family} and an isometry $V:\HH\to\KK$ such that
\begin{equation}\label{eq:extension}
 J_jV=VT_j,\qquad j\in\{1,2\}.
\end{equation}
The two-operator case of  \cite[Conjecture~5.2 (p.~275)]{BH80} is as follows.

\begin{conjecture}[Ball--Helton]\label{conj:ball-helton}
Every commuting $3$-symmetric tuple has a Jordan-family extension.
\end{conjecture}

McCullough and Russo \cite{MR16} connect positivity of a modified quadratic
symbol with a one-operator Jordan lifting by using operator Fej\'er--Riesz
factorization and completely positive maps. Their functional calculus
argument also recovers the extension theorem for a single $3$-symmetric
operator. Russo \cite{Russo17} develops the corresponding factorization
criterion for commuting tuples and uses Choi's positive map to construct a
quadratic matrix polynomial which does not factor. The construction here,
inspired by the Sobolev space connection noted above,
transfers that obstruction to the symbol of a commuting $3$-symmetric tuple.

Polynomial factors with coefficients in $\B(\HH,\KK)$ for any
target Hilbert space $\KK$ are allowed. If a quadratic polynomial $F$ satisfies
\begin{equation}\label{eq:polynomial-factor}
 F(s,t)=G(s,t)^*G(s,t),\qquad (s,t)\in\R^2,
\end{equation}
then $G$ is affine: its highest homogeneous part would otherwise contribute
a nonzero square of degree greater than two.
The extension-to-factorization implication used in \cite{BH80,Russo17}
takes the explicit form
\begin{equation}\label{eq:extension-factor}
 Q_T(s,t)=G_J(s,t)^*G_J(s,t),\qquad
 G_J(s,t)=(I+isN_1+itN_2)V.
\end{equation}
Indeed, \eqref{eq:jordan-family} gives
$e^{i(sJ_1+tJ_2)}=e^{i(sS_1+tS_2)}(I+isN_1+itN_2)$,
and $e^{i(sS_1+tS_2)}$ is unitary for real $s,t$.

For an open set $D\subset\R^2$, let $H^1(D;\C^3)$ be the space of
$L^2$ functions whose first weak derivatives belong to $L^2$, with norm
\begin{equation}\label{eq:sobolev-space-norm}
 \|f\|_{H^1(D)}^2
 =\int_D\bigl(\|f\|^2+\|\partial_1f\|^2+\|\partial_2f\|^2\bigr)\,dx.
\end{equation}
The space $H^1(D;\C^3)$ is complete in
\eqref{eq:sobolev-space-norm} \cite[Theorem 2, pp. 262–263]{Evans}, so its associated inner product makes it a
Hilbert space. Let $H_0^1(D;\C^3)$ be the closure of
$C_c^\infty(D;\C^3)$ in this norm; it is a closed Hilbert subspace.

\begin{theorem}\label{thm:counterexample}
There is an $R>0$ such that, for $D=\{x\in \R^2: \|x\|<R\},$ there is an equivalent Hilbert space
norm on $H_0^1(D;\C^3)$ for which the coordinate multipliers
$T_jf(x)=x_jf(x)$ form a bounded commuting $3$-symmetric tuple whose symbol
has no polynomial factorization \eqref{eq:polynomial-factor} into any
Hilbert space. In particular, the tuple has no Jordan-family extension.
\end{theorem}

\noindent{\bf Use of AI.} The authors acknowledge the use of AI. GPT-6.1 generated the counterexample. The authors validated the argument and modestly improved the argument and presentation. They accept responsibility for mathematical correctness.

\section{A positive matrix quadratic}
\label{sec:matrix-quadratic}

Let $\Sigma$ denote the cone of Hermitian $3\times 3$ matrix polynomials
that  have an affine factorization
\begin{equation}\label{eq:matrix-factor}
 P(s,t)=(A_0+sA_1+tA_2)^*(A_0+sA_1+tA_2) =\sum_{0\le i+j\le 2} P_{ij}s^i t^j,
 \qquad A_j\in\C^{9\times 3}.
\end{equation}
Replacing  $\C^{9\times 3}$ by bounded operators from $\mathbb{C}^3$
into a Hilbert space of dimension larger than $9$
in \eqref{eq:matrix-factor} does not enlarge $\Sigma$: the ranges of $A_0, A_1, A_2$
span a space of dimension at most $9$.

We use the standard finite-dimensional Gram representation
\begin{equation}\label{eq:gram}
 P(s,t)=Z(s,t)^*MZ(s,t),\qquad
 Z(s,t)=\begin{pmatrix}I_3\\sI_3\\tI_3\end{pmatrix}, \qquad M\succeq0.
\end{equation}
Note, from equation~\eqref{eq:matrix-factor}, one choice for $M$ is 
\[
  0\preceq M = \begin{pmatrix} A_0^*\\A_1^*\\A_2^*\end{pmatrix}
   \begin{pmatrix} A_0 & A_1 & A_2 \end{pmatrix}.
\]

The cone $\Sigma$ is closed (in any norm topology). 
Indeed, 
writing $M=[M_{ij}]_{i,j=0}^2$ in
\eqref{eq:gram}, its diagonal blocks are the coefficients $P_{00},P_{20},P_{02}$,
so $\tr M=\tr(P_{00}+P_{20}+P_{02})$. Thus representing positive Gram
matrices for a convergent sequence of quadratic polynomials have a
convergent subsequence that is again a positive Gram matrix.

Choi's example \cite{Choi75}, in the form used by Russo \cite{Russo17}, is
\begin{equation}\label{eq:choi}
 \Gamma(s,t)=
 \begin{pmatrix}
  1+2s^2&-s&-t\\
  -s&s^2+2t^2&-st\\
  -t&-st&t^2+2
 \end{pmatrix}.
\end{equation}
It is positive semidefinite on $\R^2$ and does not belong to $\Sigma$.
The nonfactorability holds over $\C$ as well as over $\R$: taking real
parts of a complex matrix sum of squares for a real symmetric polynomial
gives a real matrix sum of squares. 
\begin{lemma}\label{lem:coercive-seed}
There is a Hermitian quadratic $\Psi\notin\Sigma$ and constants $a,b>0$ such
that
\begin{equation}\label{eq:coercive-bounds}
 a(1+s^2+t^2)I_3\preceq \Psi(s,t)
 \preceq b(1+s^2+t^2)I_3,\qquad (s,t)\in\R^2.
\end{equation}
\end{lemma}

\begin{proof}
 Closedness of $\Sigma$ and \eqref{eq:choi} imply that  
\begin{equation}\label{eq:seed}
  \Gamma_\varepsilon(s,t)=\Gamma(s,t)+\varepsilon(1+s^2+t^2)I_3  \notin\Sigma
\end{equation}
for all sufficiently small $\varepsilon>0.$
Since $\Gamma\succeq0$ and $\tr \Gamma=3(1+s^2+t^2)$, \eqref{eq:coercive-bounds} holds, for such $\varepsilon,$  with $\Psi=\Gamma_\varepsilon$ 
and $a=\varepsilon$ and $b=3+\varepsilon$.
\end{proof}

\section{Sobolev realization of the quadratic}
\label{sec:realization}
Choose $\Psi$ as in Lemma~\ref{lem:coercive-seed} and write
\index{$\Psi_{ij}$}
\begin{equation}\label{eq:psi-coefficients}
 \Psi(s,t) =\sum_{0\le i+j \le 2} \Psi_{i,j} s^i t^j.
\end{equation}
Let
$D_R=\{x\in\R^2: \|x\|<R\}.$ Functions  $f\in \HH_R=H_0^1(D_R;\C^3)$ are identified
with their extension by zero to all of $\R^2$. In the usual Sobolev norm,
this embeds $H_0^1(D_R;\C^3)$ isometrically into $H^1(\R^2;\C^3)$.
Equip $\HH_R$ with the norm $\|\cdot\|_{\Psi,R}$ arising from the inner product
\begin{equation}\label{eq:sobolev-norm}
\langle f, g\rangle_{\Psi,R} = \int_{\R^2}\widehat g(\xi)^*\Psi(\xi)\widehat f(\xi)\,d\xi.
\end{equation}
Here $\widehat f$ denotes the unitary Fourier transform
\begin{equation}\label{eq:fourier-transform}
 \widehat f(\xi)=\frac{1}{2\pi}\int_{\R^2}e^{-ix\cdot\xi}f(x)\,dx
\end{equation}
of the zero extension of $f$. 
Plancherel's theorem, equations \eqref{eq:coercive-bounds} and \eqref{eq:sobolev-space-norm}  and the fact that $\widehat{\partial_jf}(\xi)=i\xi_j\widehat f(\xi)$ give
\begin{equation}\label{eq:norm-equivalence}
 a\|f\|_{H^1(\R^2)}^2\le\|f\|_{\Psi,R}^2
 \le b\|f\|_{H^1(\R^2)}^2.
\end{equation}
Consequently $\HH_R$ is complete in $\|\cdot\|_{\Psi,R}$.

Set $T_{j,R}f(x)=x_jf(x)$. The Leibniz rule and
\eqref{eq:norm-equivalence} show that $T_{j,R}$ is bounded on $\HH_R$
with norm at most $\sqrt{\frac ba}R.$ (It is easy to verify the inequality $\|T_{j,R}\|\leq \sqrt{\frac{2b}{a}} (R+1)$. A more refined argument yields the improved estimate.)
The operators $T_{j,R}$ commute. Let 
\begin{equation}\label{eq:multiplication-group}
 U_R(z,w)=e^{i(zT_{1,R}+wT_{2,R})},\qquad
 (U_R(z,w)f)(x)=e^{i(zx_1+wx_2)}f(x),
\end{equation}
where $(z,w) \in \C^2.$
For real $s,t$, write $F_R(s,t)=U_R(s,t)^*U_R(s,t)$. Thus,
letting $T_R=(T_{1,R},T_{2,R}),$ \index{$T_R$} \index{$F(s,t)$}
\begin{equation}
\label{eq:F-is-symbol}
 F_R(s,t) = Q_{T_R}(s,t)
\end{equation} 
is the symbol of $T_R$ in the sense of \eqref{eq:symbol}
and the adjoints on $\HH_R$ refer to the inner product determined by
\eqref{eq:sobolev-norm}.

\begin{proposition}\label{prop:realization}
For every $R>0$, the pair $T_R$ is a commuting $3$-symmetric
tuple. For all sufficiently large $R$, the polynomial $F_R$ has no
polynomial factorization into any Hilbert space.
\end{proposition}

\begin{proof}
Multiplication by $e^{i(sx_1+tx_2)}$ translates the Fourier transform: letting
$h=U_R(s,t)f,$ 
\begin{equation}\label{eq:fourier-translation}
 \widehat{h}(\xi)=\widehat f\bigl(\xi-(s,t)\bigr).
\end{equation}
Substituting \eqref{eq:fourier-translation} into \eqref{eq:sobolev-norm}
and setting $\eta=\xi-(s,t)$ gives
\begin{equation}\label{eq:translated-energy}
 \langle U_R(s,t)f, U_R(s,t)g\rangle_{\Psi,R}
 =\int_{\R^2}\widehat g(\xi)^*
       \Psi\bigl(\xi+(s,t)\bigr)\widehat f(\xi)\,d\xi.
\end{equation}
The expansion \eqref{eq:psi-coefficients} gives
\begin{equation}\label{eq:translated-polynomial}
\begin{aligned}
 \Psi\bigl(\xi+(s,t)\bigr)
 ={}&\Psi(\xi)
 +s\bigl(\Psi_{1,0}+2\xi_1\Psi_{2,0}+\xi_2\Psi_{1,1}\bigr)\\
 &+t\bigl(\Psi_{0,1}+\xi_1\Psi_{1,1}+2\xi_2\Psi_{0,2}\bigr)\\
 &+s^2\Psi_{2,0}+st\Psi_{1,1}+t^2\Psi_{0,2}.
\end{aligned}
\end{equation}
Every matrix coefficient $C(\xi)$ in \eqref{eq:translated-polynomial}
satisfies $\|C(\xi)\|\le c_C(1+|\xi|^2)$. Cauchy--Schwarz and
\eqref{eq:norm-equivalence} imply
\begin{equation}\label{eq:coefficient-form-bound}
 \left|\int_{\R^2}\widehat g(\xi)^*C(\xi)\widehat f(\xi)\,d\xi\right|
 \le\frac{c_C}{a}\|f\|_{\Psi,R}\|g\|_{\Psi,R}.
\end{equation}
Thus each coefficient form in \eqref{eq:translated-energy} is represented
by a bounded operator on $\HH_R$. Since $F_R=U_R^*U_R$,
\eqref{eq:translated-energy} and \eqref{eq:translated-polynomial} make
$F_R$ a $\B(\HH_R)$-valued polynomial of total degree at most two.
By \eqref{eq:F-is-symbol}, $T_R$ is a commuting $3$-symmetric tuple.

Choose a real function $\phi\in C_c^\infty(D_1)$ with $\|\phi\|_2=1$, and
define
\begin{equation}\label{eq:compression-map}
 \phi_R(x)=R^{-1}\phi(x/R),\qquad
 E_R:\C^3\longrightarrow\HH_R,\qquad E_Ru=\phi_Ru.
\end{equation}
Each $\phi_R$ has $L^2$ norm one. 
Write $\Psi_{2,0},\Psi_{1,1},\Psi_{0,2}$ for the coefficients of $s^2,st,t^2$ in $\Psi$.
Because $\phi$ is real, $|\widehat\phi(\xi)|^2$ is even under
$\xi\mapsto-\xi$, so its first moments vanish. By
\eqref{eq:compression-map}, $\widehat{\phi_R}(\xi)=R\widehat\phi(R\xi)$.
Substituting $f=\phi_Ru$ and $g=\phi_Rv$ into
\eqref{eq:translated-energy} and setting $\eta=R\xi$ gives the matrix identity
\begin{equation}\label{eq:scaled-compression}
 E_R^*F_R(s,t)E_R
 =\int_{\R^2}|\widehat\phi(\eta)|^2
       \Psi\bigl((s,t)+R^{-1}\eta\bigr)\,d\eta.
\end{equation}
Expanding $\Psi$ in \eqref{eq:scaled-compression}, using
$\int_{\R^2}|\widehat\phi(\eta)|^2\,d\eta=\|\phi\|_2^2=1$ and the
vanishing first moments of $|\widehat\phi|^2$, gives the exact compression formula
\begin{equation}\label{eq:compression}
 E_R^*F_R(s,t)E_R=\Psi(s,t)+R^{-2}K,
\end{equation}
where the fixed Hermitian matrix $K$ is
\begin{equation}\label{eq:compression-error}
 K=\int_{\R^2}|\widehat\phi(\xi)|^2
       \bigl(\xi_1^2\Psi_{2,0}+\xi_1\xi_2\Psi_{1,1}+\xi_2^2\Psi_{0,2}\bigr)\,d\xi.
\end{equation}
By \eqref{eq:seed} and closedness of $\Sigma$, the polynomial on the
right side of \eqref{eq:compression} is outside $\Sigma$ for every
sufficiently large $R$.

If $F_R=G^*G$ for an operator polynomial $G$ with any target Hilbert space,
then $G$ is affine and
\begin{equation}\label{eq:compressed-factor}
 E_R^*F_RE_R=(GE_R)^*(GE_R).
\end{equation}
The ranges of the three coefficients of $GE_R$ span a finite-dimensional
space, so \eqref{eq:compressed-factor} places
$\Psi+R^{-2}K$ in $\Sigma$, contradicting \eqref{eq:compression}.
\end{proof}

\begin{proof}[Proof of Theorem~\ref{thm:counterexample}]
Choose $\Psi$ from Lemma~\ref{lem:coercive-seed} and fix $R$ sufficiently
large for Proposition~\ref{prop:realization}. The norm induced by
\eqref{eq:sobolev-norm} and the coordinate multipliers have the stated
properties. A Jordan-family extension would give the forbidden affine
factorization \eqref{eq:extension-factor}. Hence
Conjecture~\ref{conj:ball-helton} fails.
\end{proof}

\section{The modified symbol}
\label{sec:modified-symbol}

Theorem 9 in \cite{Russo17} implicitly suggests a stronger form of the Ball-Helton conjecture. This stronger conjecture is also false as a consequence of Theorem \ref{thm:counterexample}. We explain this now.

Fix $R$ as in Proposition~\ref{prop:realization}, and write
\begin{equation}\label{eq:symbol-coefficients}
 F_R(s,t)=I+sB_{10}+tB_{01}+s^2B_{20}+stB_{11}+t^2B_{02}.
\end{equation}
Russo's lifting criterion \cite{Russo17}, extending the factorization
method of McCullough--Russo \cite{MR16}, uses
\begin{equation}\label{eq:modified-symbol}
   F_{R,c,d}(s,t)
 =F_R(s,t)-c^{-2}B_{20}-d^{-2}B_{02},\qquad c,d>0.
\end{equation}
For commuting $3$-symmetric tuples, the criterion identifies affine
factorization of \eqref{eq:modified-symbol} with extension to the
corresponding canonical Jordan family satisfying the side constraints $\|N_1\|\le c$ and $\|N_2\|\le d$ as well as $N_1N_2^*=0,$  prompting the following question: does positivity of the modified symbol suffice for such an affine factorization? The following corollary gives a negative answer. 

\begin{corollary}\label{cor:modified-symbol}
The polynomial $  F_{R,c,d}$ has no polynomial factorization into
any Hilbert space for any $c,d>0$. For all sufficiently large $c,d$, it is
positive semidefinite on $\R^2$.
\end{corollary}

\begin{proof}
The covariance argument of \cite{MR16,Russo17} gives
\begin{equation}\label{eq:modified-covariance}
   F_{R,c,d}(s,t)
 =U_R(s,t)^*\bigl(I-c^{-2}B_{20}-d^{-2}B_{02}\bigr)U_R(s,t).
\end{equation}
Choose $c,d$ so that $c^{-2}\|B_{20}\|+d^{-2}\|B_{02}\|<1$;
positivity follows from \eqref{eq:modified-covariance}.
Both $B_{20}$ and $B_{02}$ are positive. Thus a factorization of
\eqref{eq:modified-symbol}, together with constant factors for
$c^{-2}B_{20}$ and $d^{-2}B_{02}$, would factor $F_R$, contrary to
Proposition~\ref{prop:realization}.
\end{proof}

\end{document}